\documentclass[runningheads]{llncs}
\usepackage[T1]{fontenc}
\usepackage{graphicx}
\usepackage{hyperref}
\usepackage{mathtools}
\usepackage{cleveref}
\usepackage{dsfont}
\usepackage[table]{xcolor}
\usepackage{bm}
\usepackage{amsfonts}

\newcommand{\td}{\dot t}
\newcommand{\tp}{t^\prime}
\newcommand{\xd}{\dot x}
\newcommand{\xp}{x^\prime}
\newcommand{\pd}{\dot \phi}
\newcommand{\pp}{\phi^\prime}

\begin{document}
\title{From Exact Space-Time Symmetry Conservation to Automatic Mesh Refinement in Discrete Initial Boundary Value Problems}
\titlerunning{Space-time Symmetries and Mesh Refinement in Discrete IBVPs}
%
\author{Alexander Rothkopf\inst{1}\orcidID{0000-0002-5526-0809} \and
W.~A.~Horowitz\inst{2,3,4}\orcidID{0000-0001-5804-2738} \and\\
Jan Nordstr\"om\inst{5,6}\orcidID{0000-0002-7972-6183}}
\authorrunning{A. Rothkopf et al.}
%
\institute{Department of Physics, Korea University, Seoul 02841, Republic of Korea, \email{akrothkopf@korea.ac.kr} \and
Department of Physics, University of Cape Town, Private Bag X3, Rondebosch 7701, South Africa \and
Department of Physics, New Mexico State University, Las Cruces, New Mexico, 88003, USA \and
Theoretical Sciences Visiting Program, Okinawa Institute of Science and Technology Graduate University, Onna, 904-0495, Japan \and
Department of Mathematics, Link{\"o}ping University, SE-581 83 Link{\"o}ping, Sweden \and
Department of Mathematics and Applied Mathematics, University of Johannesburg, P.O. Box 524, Auckland Park 2006, Johannesburg, South Africa}
\maketitle              
\begin{abstract}
In this contribution we present recent developments in the formulation and solution of Initial Boundary Value Problems (IBVPs). Building upon a modern variational action formulation of classical dynamics, we treat Initial Boundary Value Problems directly on the action level, bypassing governing equations. We show that by including coordinate maps as dynamical degrees of freedom together with propagating fields two key results emerge. Space-time symmetries remain protected even after discretization, leading to an exact conservation of Noether charges even for discrete IBVPs. The dynamical nature of the coordinate maps leads to an adjustment of space-time resolution, guided by Noether charge conservation, realizing a form of automatic adaptive mesh refinement. We stress that as long as SBP operators are used for the discretization, our results are independent of whether the dynamics are solved on the action or governing equation level and hold in particular also at high order. As proof-of-principle for our approach we present its application to scalar wave-propagation in 1+1 dimensions. 
\end{abstract}

\section{Motivation}

\begin{figure}
\centering
\includegraphics[scale=0.4]{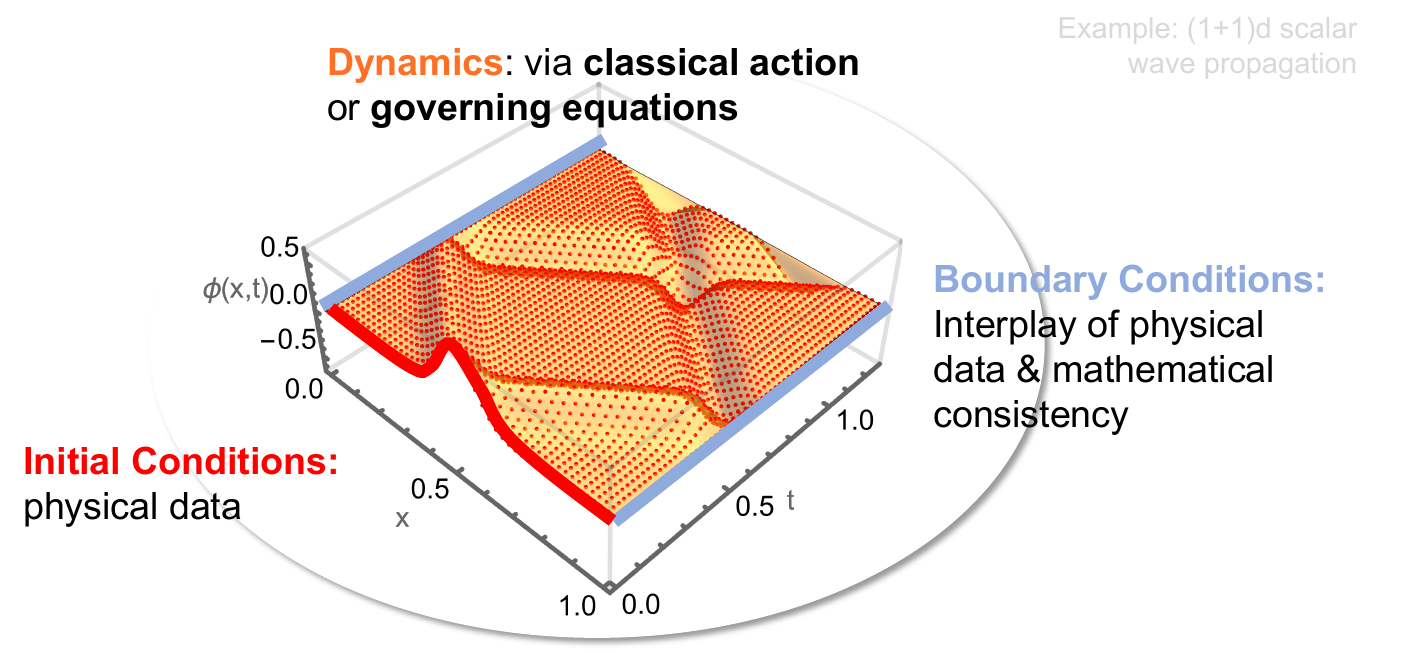}
\caption{Initial Boundary Value Problems: unique predictions for the outcome of experiments from the interplay of dynamics, informed by experimental input in the form of initial and boundary conditions.}
\end{figure}

Describing the outcome of an experiment requires us to solve an Initial Boundary Value Problem (IBVP), where we set out to predict the dynamics of a continuous degree of freedom (a field) in space and time. In the laboratory we prepare the experimental setup, representing the initial conditions for our prediction. The dynamics of the subsequent evolution in time is encoded in the classical action of the system, often represented in the form of governing equations \cite{goldstein2011classical}. These equations provide a unique prediction, if in addition to initial conditions also boundary conditions are supplied. These arise from an interplay of available physical data and mathematical consistency.

Our work is motivated by two sets of challenges that affect the formulation of IBVPs. When formulated in terms of governing equations many relevant physical systems require the use of second derivatives in time, the wave equation being one example. It is known \cite{nordstrom2016summation} that the discretization of second order systems presents a theoretical challenge in that first and second derivatives require separate discretization for a consistent treatment. Our working hypothesis is that by using an action based formulation instead, the system can be discretized using a single first order derivative. In the process of deriving governing equations one may encounter further conceptual challenges, such as the need to choose a gauge (as e.g. in electromagnetism formulated in terms of gauge potentials) and the possibility to lose intrinsic constraints which must be reintroduced a posteriori (see e.g. \cite{pinto2016handling}).

To address these challenges we set out to \textit{formulate} the IBVP \textit{on the action level} and obtain the classical field configuration directly from the action without utilizing governing equations.

The second set of challenges is related to space-time symmetries. The discretization on the level of space-time coordinates $t,{\bf x}$ breaks the continuum symmetries (see also discussion in \cite{anerot2020noether}). As Noether's theorem relates continuous symmetries and conserved charges (such as energy- or momentum density), the breaking of symmetries implies a lack of conservation. This issue affects even symplectic schemes, which preserve energy on average but not exactly at each time step.

To address these challenges, we borrow concepts from the general theory of relativity to \textit{disentangle discretization from space-time symmetries}.

\section{Part I: Action based formulation of classical dynamics}
\label{sec:part1}
\begin{figure}
\centering
\includegraphics[scale=0.45]{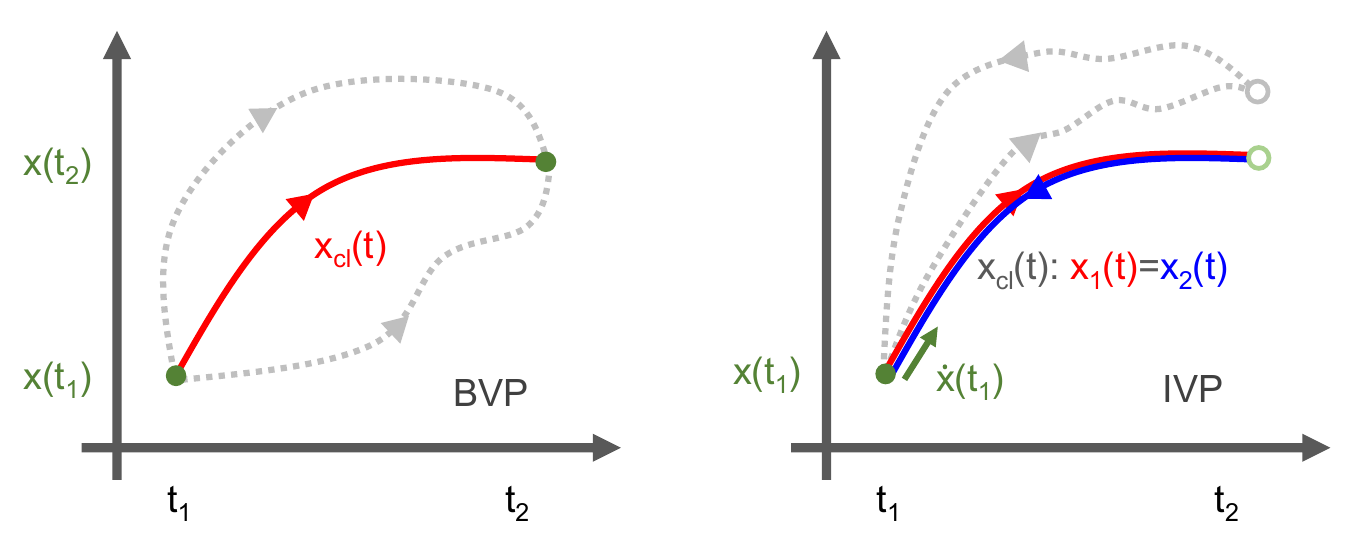}
\caption{The variational formulation of classical mechanics. (left) Traditional Hamilton's principle formulated as boundary value problem. It applies to problems where the final state of the system is known priori. (right) Modern Schwinger-Keldysh-Galley principle formulated as genuine initial value problem. It is based on a doubling of d.o.f. akin to a double shooting method. (figure from \cite{Rothkopf:2022zfb})}\label{fig:varprinc}
\end{figure}

Let us start by considering the conventional variational principle of classical mechanics, Hamilton's principle. It relies on the fact that there exists a score functional  $S[x,\dot x]$, built from the Lagrangian $L[x,\dot x]$, whose extremum is equivalent to the solution of the Euler-Lagrange equations
\begin{align}
S=\int_{t_i}^{t_f} L[x,\dot x] dt = \int _{t_i}^{t_f}\left( \frac{1}{2}m \dot x^2(t)-V(x)\right)dt, \quad \delta S=0 \Leftrightarrow \frac{\partial L}{\partial x} - \partial_t\frac{\partial L}{\partial \dot x}=0.
\end{align}
This equivalence only holds if we impose initial and final time data to combat boundary terms that arise during variation. I.e. Hamilton's principle determines the trajectory between a known starting point $x(t_i)$ and end point $x(t_f)$ (see left panel of \cref{fig:varprinc}). Note that the imposition of final time data is acausal. If we wish to predict the trajectory of a point particle, its final position in principle cannot be known a priori, otherwise there would not be a need to solve for the trajectory in the first place. The reason the boundary value nature of Hamilton's principle remains hidden in practice is that the Euler-Lagrange equations derived from Hamilton's principle are amenable to an evaluation as initial value problem. Strictly speaking, at this point, we are however only allowed to solve the governing equations as boundary value problem.

We may ask whether it is possible to formulate a variational principle on the action level for initial value problems. The problem in Hamilton's principle manifests as presence of acausal boundary terms. Since in an IVP we do not know the final position $x(t_2)$ it cannot be referenced. The solution to this conundrum has long been known in quantum theory and has recently been rediscovered in the classical context \cite{galley_classical_2013}. The Schwinger-Keldysh-Galley (SKG) action principle (for a detailed derivation see \cite{Horowitz:2026pbz}) introduces a copy of the physical degrees of freedom, placed on a forward and backward time path. This entails that the action for the doubled degrees of freedom contains the Lagrangian on the forward path and minus the Lagrangian on the backward path 
\begin{align}
S_{\rm SKG}=\int_{t_i}^{t_f} L[{\bf x}_1,\dot {\bf x}_1] - L[{\bf x}_2,\dot {\bf x}_2] dt.
\end{align}
The equations of motion are most intuitively obtained after going over to a new set of coordinates, the following linear combinations of the doubled d.o.f.
\begin{align}
{\bf x}_+=({\bf x}_1+{\bf x}_2)/2, \qquad {\bf x}_-=({\bf x}_1-{\bf x}_2)
\end{align}
The classical equations of motion follow again from the variation of this action. To see how this resembles a double-shooting method consider the right panel of \cref{fig:varprinc}. By supplying \textit{initial conditions} at $t_i$ we perform a forward shoot leading to some final value at $t_f$. At this point we use the final position and velocity as initial conditions for the backward shoot. It turns out that the classical trajectory is found if both forward and backward shoot agree. In order to initialize the backward shoot, we must identify the value and velocity of both paths at final time, leading us to the so-called \textit{connecting conditions} at $t_f$. Note that we have not fixed that value of the classical path at final time, avoiding the need to supply acausal data. Expressed in terms of the $+/-$ coordinates we have
\begin{align}
\underbracket{{\bf x}_+(t_i)=q_{\rm cl}(t_i),\; \dot {\bf x}_+(t_i)=\dot q_{\rm cl}(t_i)}_{\rm initial\, conditions}, \quad \underbracket{{\bf x}_-(t_f)=0, \;\dot {\bf x}_+(t_f)=0}_{\rm connecting \, conditions}
\end{align}
The equations of motion for ${\bf x}_+$ follow from the variation of the action w.r.t. to ${\bf x}_-$ and vice versa. What one finds is that the governing equations for ${\bf x}_-$ always allow for the vanishing constant solution, which is enforced by the connecting conditions. I.e. any term in the equation of motion for ${\bf x}_+$ that contains reference to ${\bf x}_-$ will vanish. In turn the equations of motion for the classical trajectory ${\bf x}_{\rm cl}$ can be directly obtained in the SKG formalism by the stationarity condition
\begin{align}
\left.\frac{\delta S}{\delta {\bf x}_-}\right|_{ {\bf x}_-=0, {\bf x}_+={\bf x}_{\rm cl}}=0
\end{align}
This prescription in the SKG action principle recovers the standard Euler-Lagrange equations for ${\bf x}_+$ while requiring specification of only physical initial data.

As the next step we wish to discretize this action to make it amenable to numerical optimization. This will allow us to obtain the classical trajectory without having to formulate the governing equations. In the following our focus will lie on symmetries and Noether's theorem, whose derivation in the continuum relies heavily on integration by parts (IBP). Thus we wish to retain a discrete form of IBP, which leads us to the adoption of Summation-by-Parts (SBP) finite difference operators (for reviews see e.g. \cite{svard2014review,fernandez2014review,lundquist2014sbp}). IBP connects integration and differentiation in the presence of boundaries so that in the discrete setting quadrature and finite difference stencil structure must be considered together. Let us briefly review the key ingredients to the SBP construction. We start by selecting a quadrature for discretizing the integral on $N$ grid points, encoded in a diagonal positive definite matrix $\mathds{H}$
\begin{align}
&\int dt u(t)v(t) \approx {\bf u}^t \,{\mathds H} \,{\bf v}, \quad \mathds{D}=\mathds{H}^{-1}\,\mathds{Q},\label{eq:SBPsetup}\\
&\mathds{Q}+\mathds{Q}^t=\mathds{E}_N-\mathds{E}_1={\rm diag}[-1,0,\ldots,0,1]\label{eq:SBPcond}.
\end{align}
The inverse of $\mathds{H}$ is used to define the finite difference operator $\mathds{D}$ using a stencil matrix $\mathds{Q}$, which fulfills the SBP condition given in \cref{eq:SBPcond}. In fact such an operator produces not only the necessary change in sign but also the correct boundary terms, when being moved from one side of an inner product to the other side
\begin{align}
\big({\mathds{D} \bf u}\big)^t \,{\mathds H} \,{\bf v} = - {\bf u}^t \,{\mathds H} \,{ \mathds{D}\bf v} + {\bf u}_N{\bf v}_N - {\bf u}_0{\bf v}_0.
\end{align}
The lowest order realization of the SBP principle is obtained by choosing the trapezoid rule as quadrature
\begin{align}
\mathds{H}^{[1,2,1]}=\Delta t\left[\begin{array}{cccc} \frac{1}{2}&&&\\ &1&&\\&&1&\\&&&\frac{1}{2}\end{array}\right], \quad \mathds{D}^{[1,2,1]}=\frac{1}{\Delta t}\left[\begin{array}{cccc} -1&1&0&0\\ -\frac{1}{2}&0&\frac{1}{2}&0\\0&-\frac{1}{2}&0&\frac{1}{2}\\0&0&-1&1\end{array}\right].
\end{align}

While the SBP operator defined as such is null-space consistent when considered in the context of governing equations (see \cite{svard2019convergence}), we have to consider its null-space structure more carefully if it is to be used for the discretization of an action. Let us demonstrate the issue by considering the simplest non-trivial mechanical system of a point mass in a constant force field with action $S=\int dt ( \frac{1}{2} m \dot x^2 - mg x)$. The corresponding discretized IVP action reads
\begin{align}
\nonumber \mathds{S}_{\rm IVP}&= \underbracket{\Big\{  \frac{1}{2} (\mathds{D}{\bf x}_1)^{\rm T} \mathds{H} (\mathds{D}{\bf x}_1) - g \mathds{1}^{\rm T} \mathds{H} {\bf x}_1\Big\}}_{\rm forward\,branch} - \underbracket{\Big\{\frac{1}{2} (\mathds{D}{\bf x}_2)^{\rm T} \mathds{H} (\mathds{D}{\bf x}_2) - g \mathds{1}^{\rm T}  \mathds{H} {\bf x}_2 \Big\}}_{\rm backward\,branch}\\
 \nonumber &+ \lambda_1 (  (x_1(0)+x_2(0))/2 - x_i) + \lambda_2(((\mathds{D}{\bf x}_1)(0)+(\mathds{D}{\bf x}_2)(0))/2-\dot x_i) \\
&+ \lambda_3 (x_1(N_t)-x_2(N_t)) + \lambda_4 ( (\mathds{D}{\bf x}_1)(N_t)-(\mathds{D}{\bf x}_2)(N_t) ).\label{eq:IVPfunc}
\end{align}
Here we have introduced the discretized paths ${\bf x}_{1,2}$ and make explicit both initial and connecting conditions via Lagrange multipliers (for the role of multipliers see \cite{Rothkopf:2023vki}). We will carry out a concurrent optimization in $\{ {\bf x}_{1},{\bf x}_{2},\lambda_{1-4}\}$, which amounts to a fully globally implicit scheme. As check of correct convergence we test that ${\bf x}_{1}={\bf x}_{2}$.

\begin{figure}
\centering
\includegraphics[scale=0.18]{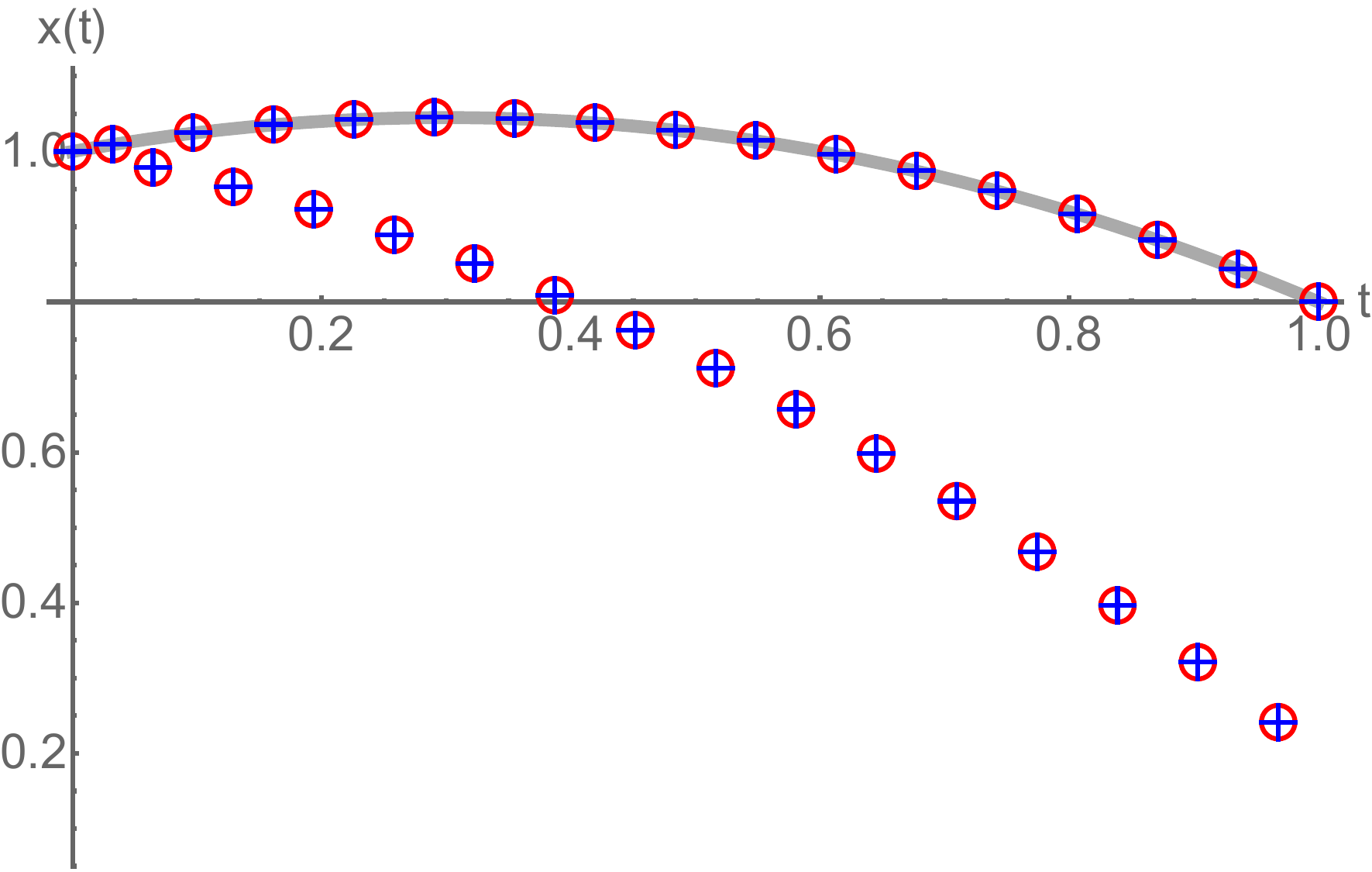}\hspace{0.7cm}
\includegraphics[scale=0.65]{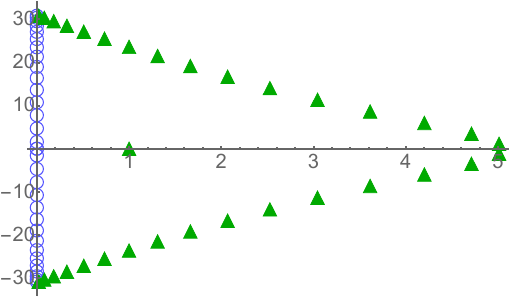}
\caption{(left) Solution for ${\bf x}_{1}$ (red open circles) and ${\bf x}_{2}$ (blue crosses), obtained from finding the extremum of the naively discretized action \cref{eq:IVPfunc} shows contamination by the highly oscillatory $\pi$-mode. Analytic solution given as gray solid line. (right) Spectrum of the unregularized lowest order SBP operator $\mathds{D}^{[1,2,1]}$ (blue open circles) and the regularized SBP operator (green triangles). Note that the latter does not show any zero mode and that the physical constant function is associated with the eigenvalue unity. (figure from \cite{Rothkopf:2022zfb})}\label{fig:solunregSBP}
\end{figure}

As shown in the right panel of \cref{fig:solunregSBP} the solution we obtain (${\bf x}_{1}$ as red open circles, ${\bf x}_{2}$ as blue crosses) only follows the analytic solution (gray solid) closely at half the time steps. The other half exhibits contamination by the so-called $\pi$-mode the maximally oscillating function supported on a grid of size N. The reason for the occurrence of the $\pi$-mode lies in the spectrum of the SBP operator, shown as blue open circles in the right panel of \cref{fig:solunregSBP}. It features two zero eigenvalues. When inspecting the right eigenvectors, both zeros are associated with the constant function, the two zero modes are degenerate. Since $\mathds{D}$ is not antisymmetric due to the boundary its left and right eigenvectors do not necessarily agree. And indeed the left eigenvectors associated with the zero modes are indeed the $\pi$-mode. They too are degenerate. 

While in a differential equation setting only the right eigenvectors are relevant, here terms like $(\mathds{D}{\bf x})^t \mathds{H} \,\mathds{D}{\bf x}$ appear in the action, where both left and right eigenvectors may play a role in determining the extremum. 

\begin{figure}
\centering
\includegraphics[scale=0.19]{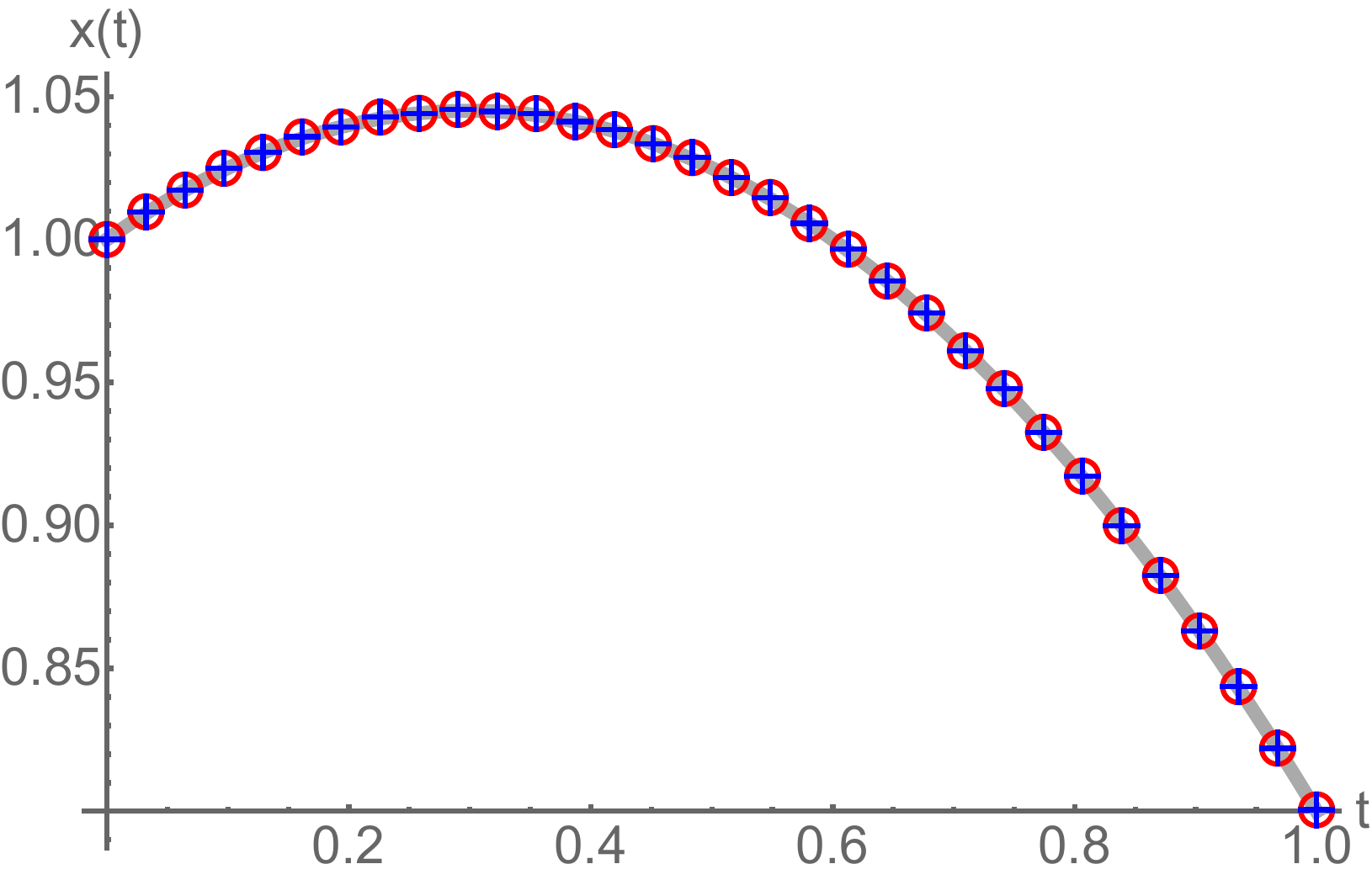}\hspace{0.7cm}
\includegraphics[scale=0.24]{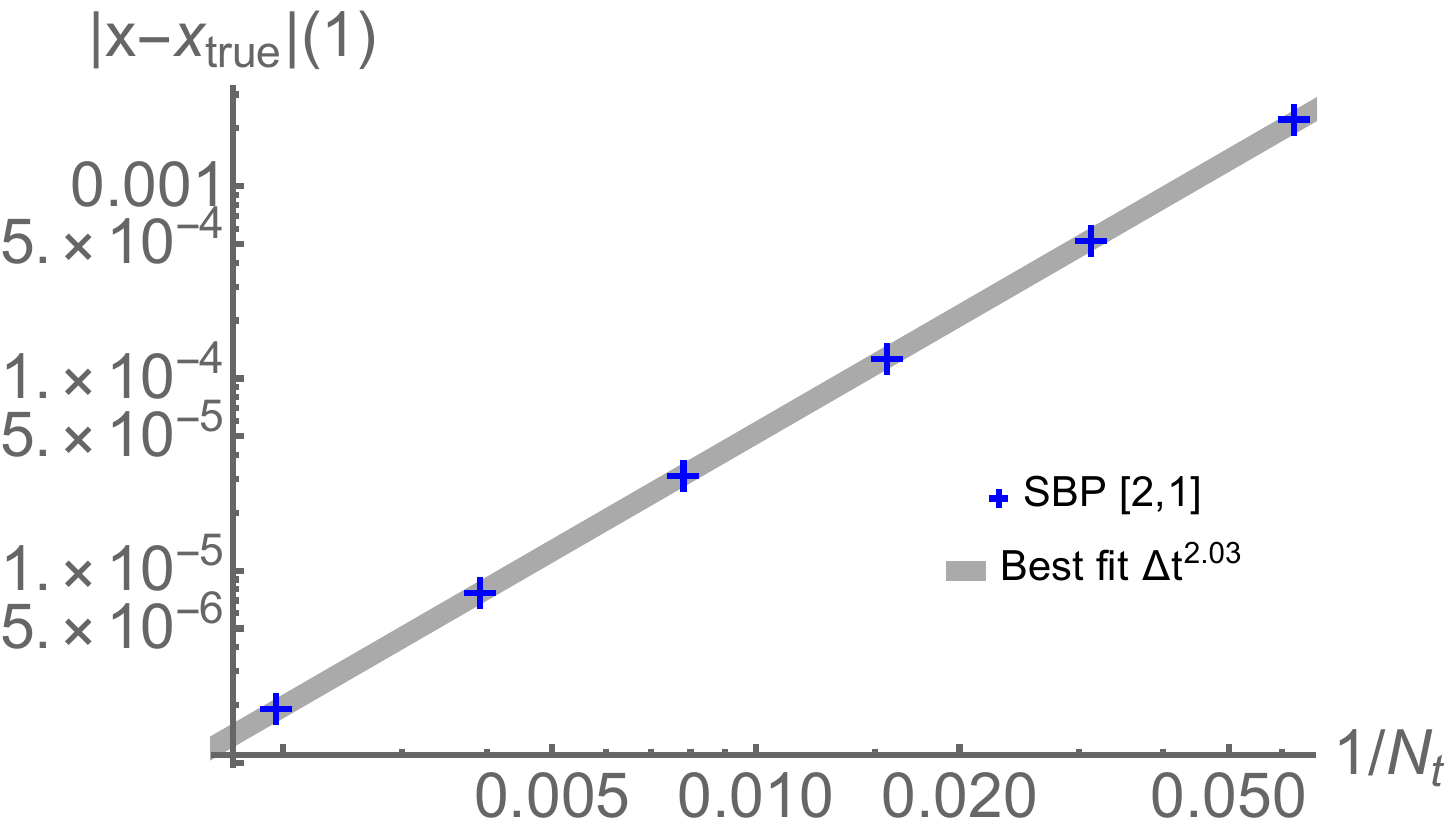}
\caption{(left) Solution for ${\bf x}_{1}$ (red open circles) and ${\bf x}_{2}$ (blue crosses), obtained from finding the extremum of the discretized action \cref{eq:IVPfunc} using the regularized SBP operator $\bar{\mathds{D}}^{[1,2,1]}$. No more contamination by the $\pi$-mode appears in the comparison to the analytic solution (gray solid). (right) Competitive convergence of the solution under grid refinement showing a powerlaw behavior with $\Delta t^{2.03}$ as expected for a well-defines SBP scheme of lowest order. (figure from \cite{Rothkopf:2022zfb})}\label{fig:solregSBP}
\end{figure}

We decide to tackle the $\pi$-mode by taking inspiration from the modern treatment of initial boundary value problems, using weakly imposed boundary data to provide regularization. The Simultaneous Approximation Term (SAT) technique \cite{lundquist2014sbp} introduces penalty terms to impose boundary conditions. Here we will modify the finite difference operator $\mathds{D}$ by adding to it a penalty term using affine coordinates. Defining a discrete array ${\bf x}_0$ which only contains the initial value in its first entry, we can write
\begin{align}
\bar{\mathds{D}}{\bf x}=\mathds{D}{\bf x}+\mathds{H}^{-1}\mathds{E}_1\big({\bf x}-{\bf x}_i\big).
\end{align}
The matrix $\mathds{E}_1$ projects out the first entry of the array its is applied to. While the part of the penalty term that is proportional to ${\bf x}$ can be absorbed straight forwardly into the new $\bar{\mathds{D}}$ the constant shift requires us to extend the $N\times N$ matrix structure of $\mathds{D}$ by a single column and row as indicated by the light blue shading below
\begin{align}
\frac{1}{\Delta t}\left[ \begin{array}{cccc>{\columncolor{blue!20}}c} -1 + 2 & 1 & 0 & 0 &-2 x_i \\ -\frac{1}{2} &0 &\frac{1}{2} & 0 & 0\\ 0 & -\frac{1}{2} &0 &\frac{1}{2} &0 \\ 0&0&-1&1 &0 \\\rowcolor{blue!20} 0&0&0&0&1\end{array}\right].
\end{align}
The newly added column carries the information about the shift. When we plot the spectrum of this regularized SBP operator, we find on the right panel of \cref{fig:solregSBP} that the green triangles show no more zero modes. All eigenvalues have positive real part and, except for one, occur in complex conjugate pairs. The eigenvalue at unity is associated with the physical constant mode. Let us stress that this \textit{regularization prescription} is \textit{independent of the order} of the SBP operator and thus applies also for higher order schemes. It offers an alternative to the conventional upwind regularization \cite{courant1952solution} or the Wilson term regularization \cite{Wilson:1974sk} applicable in case of complex valued functions.

Repeating the global optimization of \cref{eq:IVPfunc} after replacing unregularized SBP operators by their regularized counterparts finally provides us with a classical trajectory devoid of the $\pi$-mode contamination, as shown in the left panel of \cref{fig:solregSBP}. The solution for ${\bf x}_{1}$ (red open circles) and ${\bf x}_{2}$ (blue crosses) now cleanly follows closely the analytic solution given as gray solid line. Performing a grid refinement analysis we confirm that the solution obtained from the extremum of the IVP action approaches the true analytic solution with a powerlaw according to $\Delta t^{2.03}$, as shown in the right panel of \cref{fig:solregSBP}. I.e. our approach exhibits scaling competitive with the conventional lowest order SBP treatment of the governing equations. In \cite{Rothkopf:2022zfb} we confirmed competitive scaling also for the SBP operator of order $[4,2,4]$, indicating that competitive scaling persists also for higher order schemes. 

\section{Part II: Dynamical coordinate maps for space-time symmetry preservation}

In the preceding section we described how to formulate IVPs directly on the action level. On the one hand this gives us direct access to a discrete form of Noether's theorem, on the other hand it does not yet provide conserved Noether charges, since a discretization in terms of space-time coordinates (in case of an IVP, i.e. time) violated continuous space-time symmetries.

To disentangle discretization from space-time symmetries let us consider the treatment of point particle mechanics in general relativity (for an introduction see \cite{carroll2019spacetime}). We will not be concerned with dynamical gravity here, we simply wish to borrow the underlying mathematical formalism to use it to our advantage. I.e. the motion of our system will always proceed in conventional flat spacetime. 

However general relativity emphasizes that space and time must be treated on the same footing. I.e. instead of thinking of the trajectory of a particle as a function ${\bf x}(t)$, we should think of the path that the particle traces out in space-time as a one-dimensional manifold, a line. This line can be parameterized by an abstract world-line parameter $\gamma$. Along the world-line the particle carries with it both a space and a time label encoded in the coordinate maps $t(\gamma)$ and ${\bf x}(\gamma)$. Often these maps are collected in a vector $X^\mu(\gamma)=(t(\gamma), {\bf x}(\gamma))^\mu$. 

General relativity tells us that the motion of a point particle is described by the \textit{shortest path in a given space-time}, a so-called geodesic. Space-time is defined in general relativity from the so-called metric tensor $G_{\mu\nu}$, which defines distances and angle relations among different space-time points. For conventional flat spac-time $G_{\mu\nu}={\rm diag}[c^2,-1,-1,-1]$.

One may ask whether geodesics too arises from a variational principle. The answer is in the affirmative and the corresponding geodesic action takes on the following form
\begin{align}
S= \int_{\gamma_i}^{\gamma_f}\, d\gamma\, (-mc)\sqrt{ \big(G_{00} + V(X)/2mc^2\big)\frac{dX^{0}}{d\gamma} \frac{dX^{0}}{d\gamma} - G_{ii} \frac{dX^{i}}{d\gamma} \frac{dX^{i}}{d\gamma} }    \label{eq:geodact}.
\end{align}
And while this formula appears involved, a closer inspection reveals that it is but the multidimensional generalization of the arc-length formula known from elementary calculus $\ell = \int \sqrt{1+f'(x)^2}$. As shown in \cite{Rothkopf:2023ljz}, if we consider motion where the potential energy remains much lower than the rest energy of the particle $V(x)/2mc^2\ll 1$ and velocity remains well below the speed of light $dX^i/d\gamma / dt/d\gamma /c = v/c \ll 1$ we can take the non-relativisitc limit and recover the standard classical action
\begin{align}
S_{\rm nr}=&\int_{t_i}^{t_f}\,dt \Big( -mc^2   + \frac{1}{2} m \Big(\frac{d{\bf x}}{dt}\Big)^2 - V(x) \Big)\label{eq:eq6}.
\end{align}
Note the appearance of the \textit{constant term} $-mc^2$ in \cref{eq:eq6}. Had we not known about the geodesic action, we would have been unable to discover this term in the non-relativistic setting, since a constant term does not contribute to the Euler-Lagrange equations. Importantly $mc^2$ denotes an important \textit{physical scale}, which is where the motion through space and motion through time become inseparable. 

The world-line action offers an interesting new handle to discretization, in that we may discretize in the abstract parameter $\gamma$ instead of time. While it entails that the coordinate maps $t(\gamma_i)$ and ${\bf x}(\gamma_i)$ are evaluated on discrete values of $\gamma_i$, the values of the map themselves remain continuous. In \cite{Rothkopf:2023ljz} we showed that since the maps remain continuous, the world-line action remains invariant under \textit{continuous space-time symmetries} even after discretization in $\gamma$. In turn Noether's theorem provides us with conserved charges in the discrete setting.

\begin{figure}[t]
\centering
\includegraphics[scale=0.3]{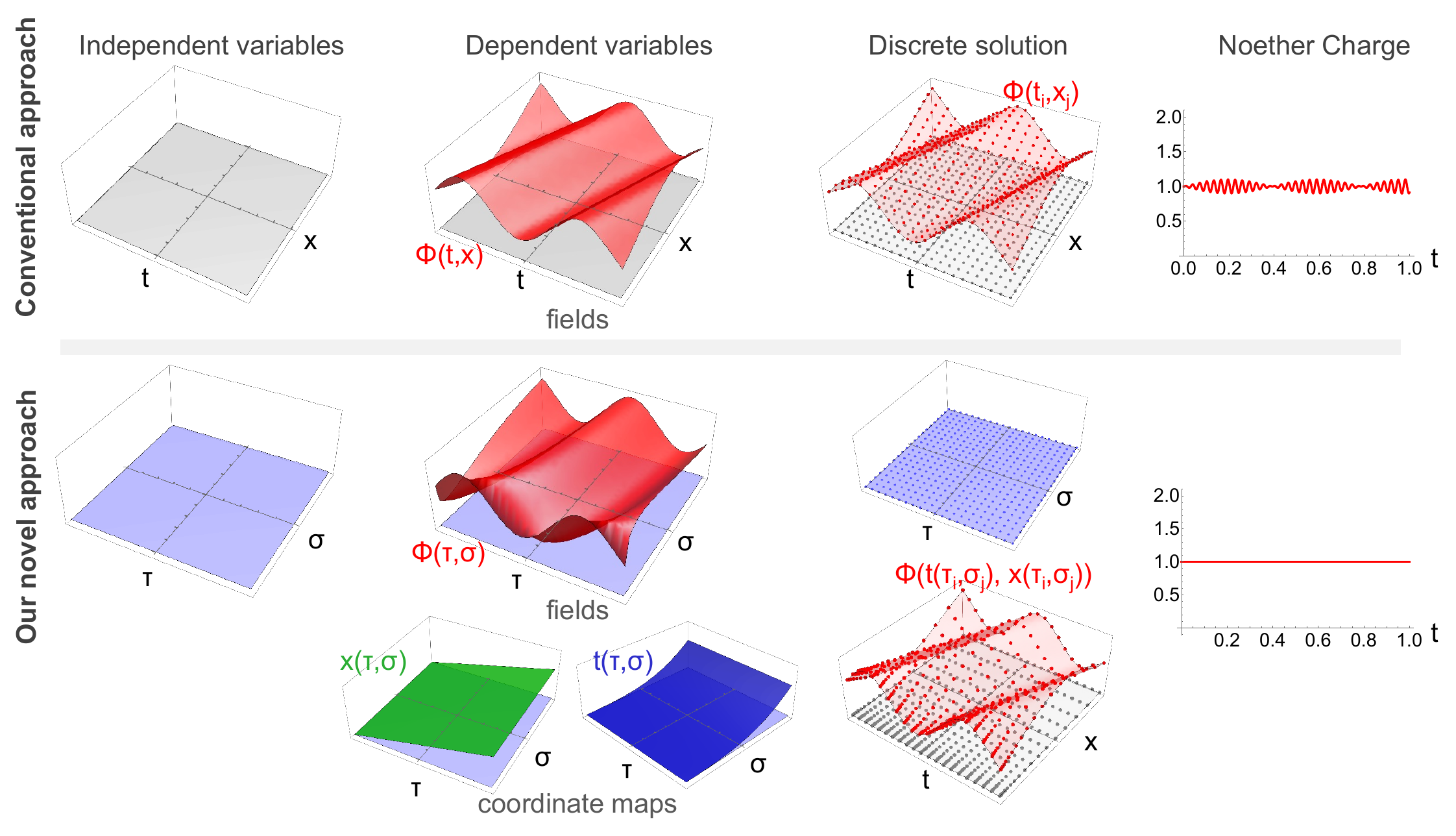}
\caption{(top row) Conventional formulation of IBVPs using space-time coordinates $(t,x)$ as independent variables. The dynamical degree of freedom $\phi(t,x)$ evolves as function of space-time coordinates. When naively discretized, the grid structure in space-time coordinates breaks space-time symmetries and continuum Noether charges fail to be conserved. (bottom row) Proposed formulation of IBVPs in the presence of dynamical coordinate maps. Inspired by the world-line formalism we consider abstract parameters $\tau$ and $\sigma$ to parameterize temporal and spatial directions. Not only the conventional dynamical d.o.f. $\phi(\tau,\sigma)$ is expressed as a function of $\tau,\sigma$ but in addition we introduce dynamical time $t(\tau,\sigma)$ and spatial $x(\tau,\sigma)$ maps that evolve alongside $\phi$. We show that discretizing $\tau,\sigma$ leaves the coordinate maps continuous and thus allows us to retain continuum space-time symmetries in the discrete setting, leading to exactly preserved Noether charges. (figure from \cite{Rothkopf:2024hxi})}\label{fig:newIBVP}
\end{figure}

In this contribution we report on the extension of this idea from IVPs to IBVPs as developed in \cite{Rothkopf:2024hxi}. The differences between the conventional formulation of IBVPs and our proposed novel approach are outlined and discussed in \cref{fig:newIBVP}.

Our starting point is the reparameterization invariant action, a standard expression in general relativity textbooks (we took the liberty to pull out an overall factor of one-half for notational convenience)
\begin{align}
    S=\int d^{(d+1)}X \, \sqrt{-{\rm det}[G]} \Big\{ \frac{1}{2}\Big(G^{\mu\nu} \partial_\mu\phi(X)\partial_\nu\phi(X) + V(\phi) \Big)\Big\}.\label{eq:stdFTaction}
\end{align}
Remembering the discussion around the non-relativistic limit of the world-line action we may ask whether this action is missing a constant term. We proceed by adding it in the form of $-T$ below
\begin{align}
    S=\int d^{(d+1)}X \, \sqrt{-{\rm det}[G]} \Big\{ -T + \frac{1}{2}\Big(G^{\mu\nu} \partial_\mu\phi(X)\partial_\nu\phi(X) + V(\phi) \Big)\Big\}.\label{eq:stdFTactionpconst}
\end{align}
So far we have not changed the classical motion. Let us factor out the term $(-T)$ up front which leads to the following expression
\begin{align}
 {\cal S}&\equiv \int d^{(d+1)}X \, \sqrt{-{\rm det}[G]} \big(-T\big) \Big\{ 1 - \frac{1}{2T}\Big(G^{\mu\nu} \partial_\mu\phi(X)\partial_\nu\phi(X) + V(\phi) \Big)\Big\} + {\cal O}(\kappa^2),
\end{align}
where the presence of an integrand of the form $\{1-\frac{1}{2} \kappa\}$ reminds us of the first terms in the series expansion we encountered in the world-line action. Our proposal amounts to interpreting this action as the low-energy limit of a more general action including higher powers of $\kappa={\rm action\, density}/T$. Undoing the apparent expansion leads us to an expression featuring a square-root term (similar to the world-line action in \cref{eq:geodact})
\begin{align}
 {\cal S}=\int d^{(d+1)}X \, \sqrt{-{\rm det}[G]} \big(-T\big) \sqrt{ 1 - \frac{1}{T}\Big(G^{\mu\nu} \partial_\mu\phi(X)\partial_\nu\phi(X) + V(\phi) \Big)}.  \label{eq:stdFTactionbacksqrt}
\end{align}
The crucial next step lies in elevating the space-time coordinates to actual dynamical coordinate maps. Whereas in the world-line action we had $t\to t(\gamma)$ we here introduce $X^\mu\to X^\mu(\tau,{\bm \sigma})$, where $\Sigma^a=(\tau,{\bm \sigma})^a$ denote the new independent variables for our integral. The change of variables $X^\mu \to \Sigma^a$ introduces the Jacobian of the coordinate transform into the integrand. The reparameterization invariance of the action however allows us to write the action in the new variables in the same form as in the old variables
\begin{align}
     {\cal S}=\int d^{(d+1)}\Sigma \, \sqrt{-{\rm det}[g]}& \big(-T\big)\sqrt{ 1 - \frac{1}{T}\Big(G^{\mu\nu} \partial_\mu\phi(X(\Sigma))\partial_\nu\phi(X(\Sigma)) + V(\phi) \Big)},\label{eq:stdFTactionbacksqrtmod1}   
\end{align}
where $\sqrt{-{\rm det}[g]}$ now makes reference to the \textit{induced metric} on the space of abstract coordinates, instead of the metric of the former space-time coordinates. Note that the induced metric is defined from the known space-time metric $G$ and the dynamical coordinate maps $g_{ab} =G_{\mu\nu} (dX^\mu/d\Sigma^a) (dX^\nu/d\Sigma^b)$.

The terms under the square root still make reference to derivatives in the old coordinates with indices $(\mu,\nu)$ but with a bit of algebra they can be re-expressed fully in terms of the new coordinates with indices $(a,b)$. The final result for our novel action in the presence of dynamical coordinate maps, introduced in \cite{Rothkopf:2024hxi}, reads
\begin{align}
{\cal S}_{\rm BVP}=\int d^{(d+1)}\Sigma \, \big(-T\big) \sqrt{ \Big(\frac{1}{T}V(\phi)-1\Big){\rm det}[g] + \frac{1}{T} \partial_a\phi(\Sigma)\partial_b\phi(\Sigma) {\rm adj}[g]_{ab}}.
\end{align}
The first term under the square root only depends on the coordinate maps and is known in the theoretical physics literature as the Nambu-Goto action \cite{zwiebach2004first}. The second term contains the field $\phi$ dependence and the adjugate of the metric ${\rm adj}[g]_{ab}$ connects the derivatives in the new coordinates to those in the old space-time coordinates. The scale $T$ in analogy to the world-line action here denotes a scale at which the dynamics of the field becomes inseparable from the dynamics of the coordinate maps.

In order to formulate the evolution of the field $\phi$ and the coordinate maps $X^\mu$ as a genuine initial value problem we must return to the Schwinger-Keldysh-Galley prescription of doubling all dynamical d.o.f. and assigning Lagrangians with opposite sign to each copy
\begin{align}
{\cal S}_{\rm IBVP}&= \int d^{(d+1)}\Sigma\; \Big\{ S_{\rm BVP}[X_1,\partial_a X_1,\phi_1, \partial_a\phi_1]-S_{\rm BVP}[X_2,\partial_a X_2,\phi_2, \partial_a\phi_2]\Big\}.
\end{align}
Here we must now supply initial and in particular connecting conditions for all propagating degrees of freedom
\begin{align}
&X_1^\mu(\tau=\tau^{\rm f},\vec{\sigma})=X_2^\mu(\tau=\tau^{\rm f},\vec{\sigma}), \quad \partial_0 X_1^\mu|_{\tau=\tau^{\rm f}}= \partial_0 X_2^\mu|_{\tau=\tau^{\rm f}},\\
&\phi_1(\tau=\tau^{\rm f},\vec{\sigma})=\phi_2(\tau=\tau^{\rm f},\vec{\sigma}), \quad  \partial_0 \phi_1|_{\tau=\tau^{\rm f}}= \partial_0 \phi_2|_{\tau=\tau^{\rm f}}.\label{eq:SKconnecting}
\end{align}
Note that if we set values equal on the final time slice, then we have also automatically set equal all of the spatial derivatives. I.e. why in the above connecting conditions only values and temporal derivatives are explicitly mentioned (for more details of the SKG construction see \cite{Rothkopf:2024hxi}).

In order to obtain the solution numerically, we must discretize the action ${\cal S}_{\rm BVP}$ and we decide to do so using the simplest generalization of the regularized one-dimensional SBP operators introduced in \cref{sec:part1}. Introducing the $(d+1)$ dimensional discretized arrays ${\bf X}^\mu$ and ${\bm \phi}$ we use Kronecker products to construct from one-dimensional SBP finite difference operators $(d+1)$ dimensional ones. While the SBP stencil structure is the same for each d.o.f. in the system, we must use appropriate penalty terms for their regularization. I.e. there exists individually regularized versions of the SBP operators to be applied to either the coordinate maps or the field, such as $\bar{\mathds{D}}_a^\mu {\bf X}^\mu$ where the repeated index is not summed over but indicates onto which spacetime coordinate the $\Sigma_a$ derivative is applied. Similarly we have $\bar{\mathds{D}}^\phi_a {\bm \phi}$ indicating the appropriately regularized derivative operator in $\Sigma_a$ direction acting on the field. 

The resulting discretized action is given as
\begin{align}
 \nonumber &\mathds{E}_{\rm BVP}[{\bm X}_1^\mu,\bar{\mathds{D}}^\mu_a {\bm X}_1^\mu,{\bm \phi}_1, \bar{\mathds{D}}^\phi_a{\bm \phi}_1]\\
 &=\frac{1}{2}\Big\{ \Big(\frac{1}{T}V({\bm \phi}_1)-1\Big)\circ{\rm det}[{\bm g}_1] +\frac{1}{T} (\bar{\mathds{D}}^\phi_a{\bm \phi}_1) \circ  (\bar{\mathds{D}}^\phi_b{\bm \phi}_1)\circ{\rm adj}[{\bm g}_1]_{ab}\Big\}^{\frac{1}{2}} {\bm h}.\label{eq:discrEBVP}
\end{align}
Here the inner product in function space has been rewritten in a simplified notation, where the discrete array $\{\ldots\}^{\frac{1}{2}}$, arising from the element-wise application of the square root, is dotted into the discrete array ${\bf h}$, which contains the diagonal entries of the quadrature matrix ${\mathds{H}}$. Let us inspect how the discrete coordinate maps enter this action. They do so via the induced metric and its adjugate
\begin{align}
     {\bm g}_{ab}=G_{\mu\nu} (\bar{\mathds{D}}^\mu_a {\bm X}^\mu)\circ (\bar{\mathds{D}}^\nu_b {\bm X}^\nu), \quad {\rm det}[{\bm g}]=\sum_{i_0,\ldots,i_{d}}\epsilon_{i_0\cdots i_{d}} {\bm g}_{0,i_0}\circ\cdots\circ{\bm g}_{d,i_{d}}\label{eq:discrg}.
\end{align}
We would like to highlight the fact that after discretization each entry of the ${\bf g}_{ab}$ matrix \textit{remains invariant under the full continuous Poincar\'e transformation group}. The fact that derivatives are contracted with the space-time metric $G_{\mu\nu}$ means that the expression is invariant under rotations and space-time boosts (the relativistic generalization of rotations involving both space and time coordinates). The fact that a derivative acts on the coordinate maps entails that constant shifts are annihilated, leading to translational invariance. Since the entries of the coordinate maps remain continuous, these symmetries can be realized infinitesimally, as required in Noether's theorem to yield conserved charges.

Since we use a mimetic SBP operator, the derivation of the Noether charge proceeds as in the continuum, simply replacing derivatives by SBP finite differences. For completeness let us quote here the explicit expression for the Noether charge associated with space-time translations, including terms that arise from the imposition of initial and connecting conditions via Lagrange multipliers $\tilde{\bm \lambda}$ and $\tilde{\bm \gamma}$. The discrete version of the one-dimensional delta function is denoted by  $\delta(t-t_i)\approx \mathfrak{d}^0[i]$
\begin{align}
    {\bm Q}^{\rm L} =  \Big( \mathds{H}_\sigma \frac{\partial {\mathds{E}}_{\rm BVP}}{\partial (\mathds{D}_0 {\bf X}^\mu)} + ({\bm h}_\sigma^T \tilde{\bm \lambda}_\mu) \mathfrak{d}^0[0] + ( {\bm h}_\sigma^T \tilde{\bm \gamma}_\mu) \mathfrak{d}^0[N_0] \Big)\delta {\bf X}^\mu.\label{eq:NoetherChargeDiscrLagr}
\end{align}

\begin{figure}[t]
\centering
\includegraphics[scale=0.475]{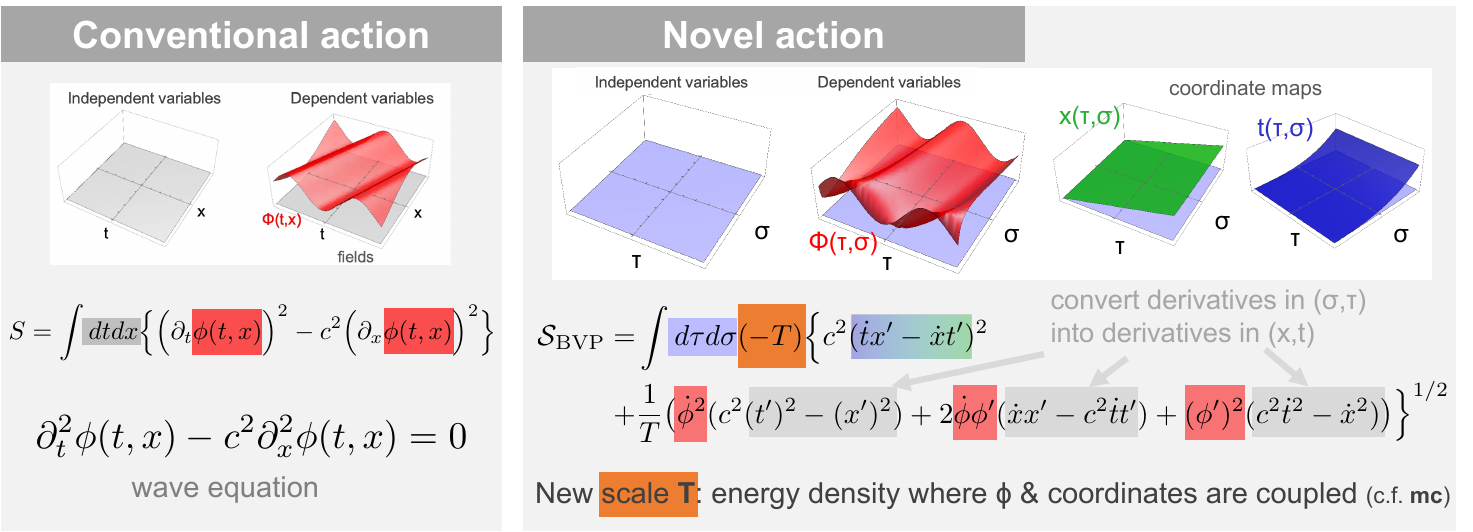}
\caption{The setup for our proof-of-principle in 1+1 dimensions. (left) The conventional action approach to scalar wave evolution. The action is expressed in terms of space-time coordinates and its variation leads to the well-known wave equation. (right) Our new action is expressed in abstract $\tau,\sigma$ variables. The Nambu-Goto term (blue-green shaded) takes on a simple form involving only coordinate map derivatives. The field dependent term contains the different derivatives of the field w.r.t. the abstract parameters, each multiplied by a term that relates $\tau,\sigma$ derivatives to the original $x,t$ derivatives. The new scale $T$ serves as a reminder that field dynamics and coordinate dynamics are now interlinked.}\label{fig:1p1dcomp}
\end{figure}

As proof-of-principle we present a treatment of scalar wave propagation in $(1+1)$ dimensions. The full action in this case reads
\begin{align}
 {\cal S}_{\rm BVP}=\int d\tau d\sigma& \, \big(-T\big) \Big\{ c^2(\td\xp-\xd\tp)^2 + \label{eq:novelaction1p1}\\
 \nonumber&\frac{1}{T} \Big( \pd^2(c^2(\tp)^2-(\xp)^2) + 2\pd\pp(\xd\xp-c^2\td\tp) + (\pp)^2(c^2\td^2-\xd^2) \Big) \Big\}^{1/2}. 
\end{align}
A comparison of this action to the conventional action of scalar wave propagation is given in \cref{fig:1p1dcomp}. In the following we will simplify the computation by restricting ourselves to only a dynamical temporal map $t=t(\tau,\sigma)$ and we fix apriori to a trivial spatial map $x(\tau,\sigma)=\sigma$.

Even though the resulting action is only invariant under infinitesimal time translations, the dynamical nature of the t-map suffices to demonstrate the validity of our approach
 \begin{align}
{\cal E}_{\rm BVP}\overset{x=\sigma}{=}\int &d\tau d\sigma \, \frac{1}{2} \Big\{ (\td)^2+\frac{1}{T} \Big( \pd^2((\tp)^2-1) - 2 \pd\pp\td\tp + (\pp)^2(\td^2) \Big) \Big\}. \label{eq:novelactionE1p1}
\end{align}
Note that we have squared the integrand here. Since taking the square amounts to a monotonic transformation, the new integrand will lead to the same extremum of the action. 

After converting the action of \cref{eq:novelactionE1p1} to its SKG form and discretizing it with regularized SBP operators on a $N_\tau\times N_\sigma = 60 \times 48$ grid, we obtain the classical solution for the field ${\bm \phi}$, as well as the dynamical coordinate maps ${\bm t}$ from a numerical optimization using the \texttt{IPOPT} library implemented in the \texttt{Mathematica} software (open access code is available at \cite{rothkopf_2024_11082746}). In the following we use $T=10.000$ in order to achieve coupling of the maps with the field dynamics on the percent level.

\begin{figure}
\centering
\includegraphics[scale=0.25]{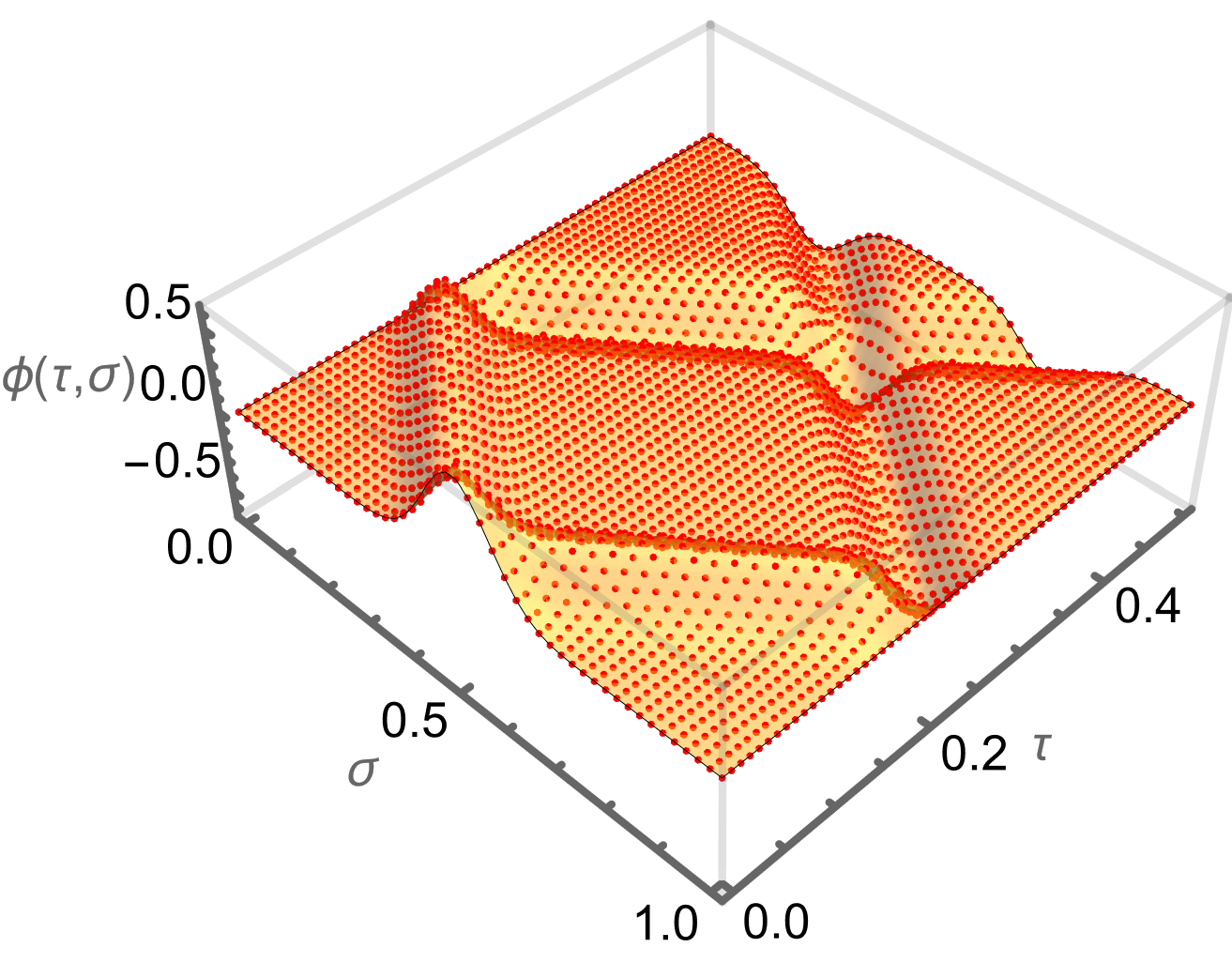}\hspace{0.7cm}
\includegraphics[scale=0.25]{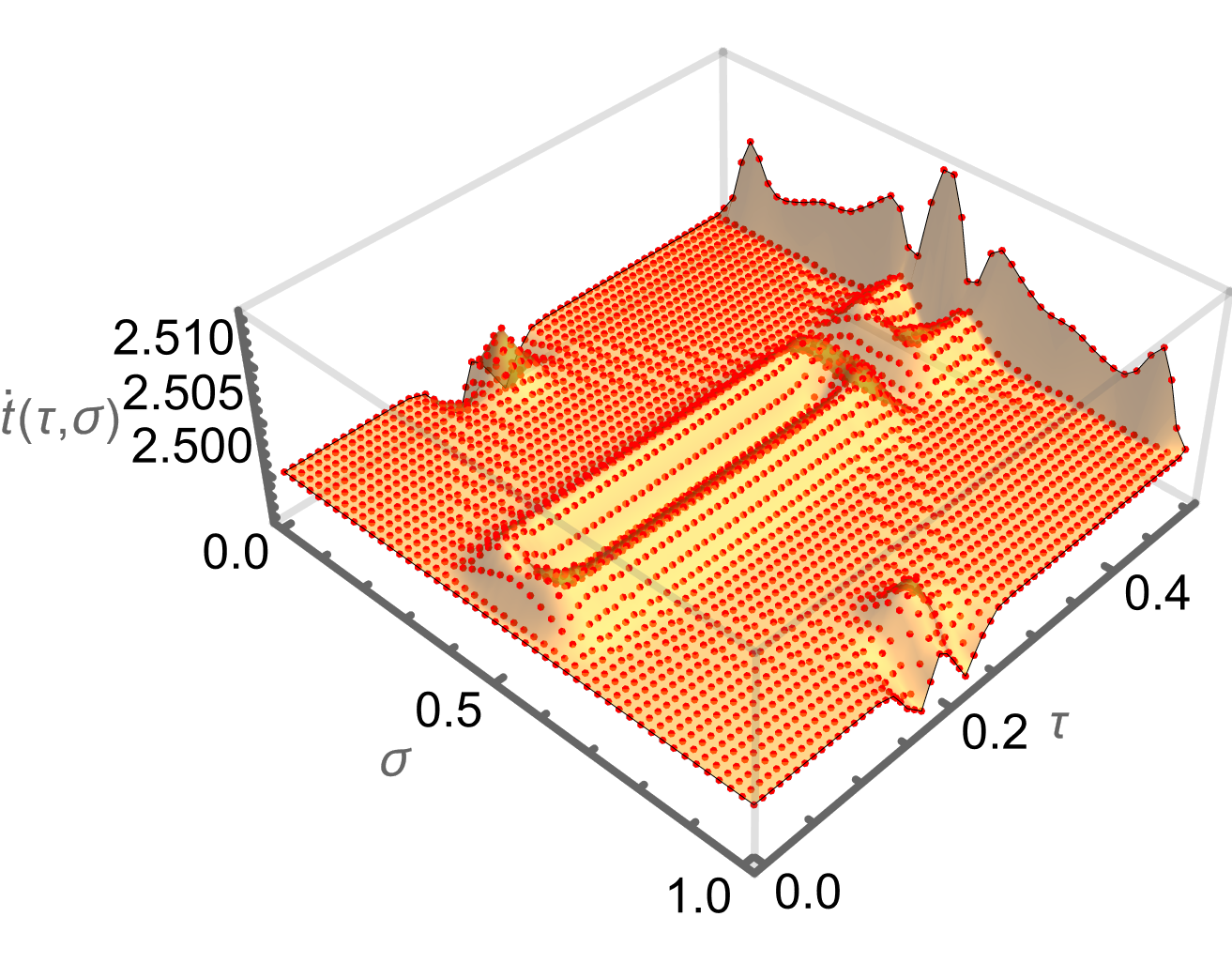}
\caption{(left) Numerical solution of the field evolution as a function of $\tau$ and $\sigma$, obtained from the optimization of the discretized action \cref{eq:novelactionE1p1}. Wave propagation is initialized from a bump centered in the spatial domain from which both a left- and right-moving wave-package emanates. The packages approach the Dirichlet boundary, reflect, including a phase flip, and approach the center of the simulation domain. When they meet, they interfere constructively and begin to move towards the opposite spatial boundary. (right) The $\tau$-derivative of the dynamical temporal map. Large values here denote a coarser time resolution, while smaller values indicate a finer time resolution. We find interesting structures that closely align with the field dynamics. As the wave packages recede from the interior of the simulation domain the dynamical coordinate map automatically reduces the resolution there. On the other hand, at the boundary, where the wave-packages violently scatter off the Dirichlet boundary, the dynamics require finer resolution as indicated by the lower values found there. (figure from \cite{Rothkopf:2024hxi})}\label{fig:numsol}
\end{figure}

As shown and described in detail in the caption of \cref{fig:numsol} we find that the dynamical coordinate map adjusts automatically to the dynamics of the propagating field. When the wave-packages in our simulation leave the interior of the simulation domain, the temporal map decreases the resolution there. On the other hand, in the instant that the wave-packages scatter from the Dirichlet boundary, we find that the temporal map introduces a finer time resolution to adequately capture the violent dynamics. This automatic adaptation of the temporal resolution amounts to a intrinsic form of \textit{automatic adaptive mesh refinement}.

\begin{figure}
\centering
\includegraphics[scale=0.25]{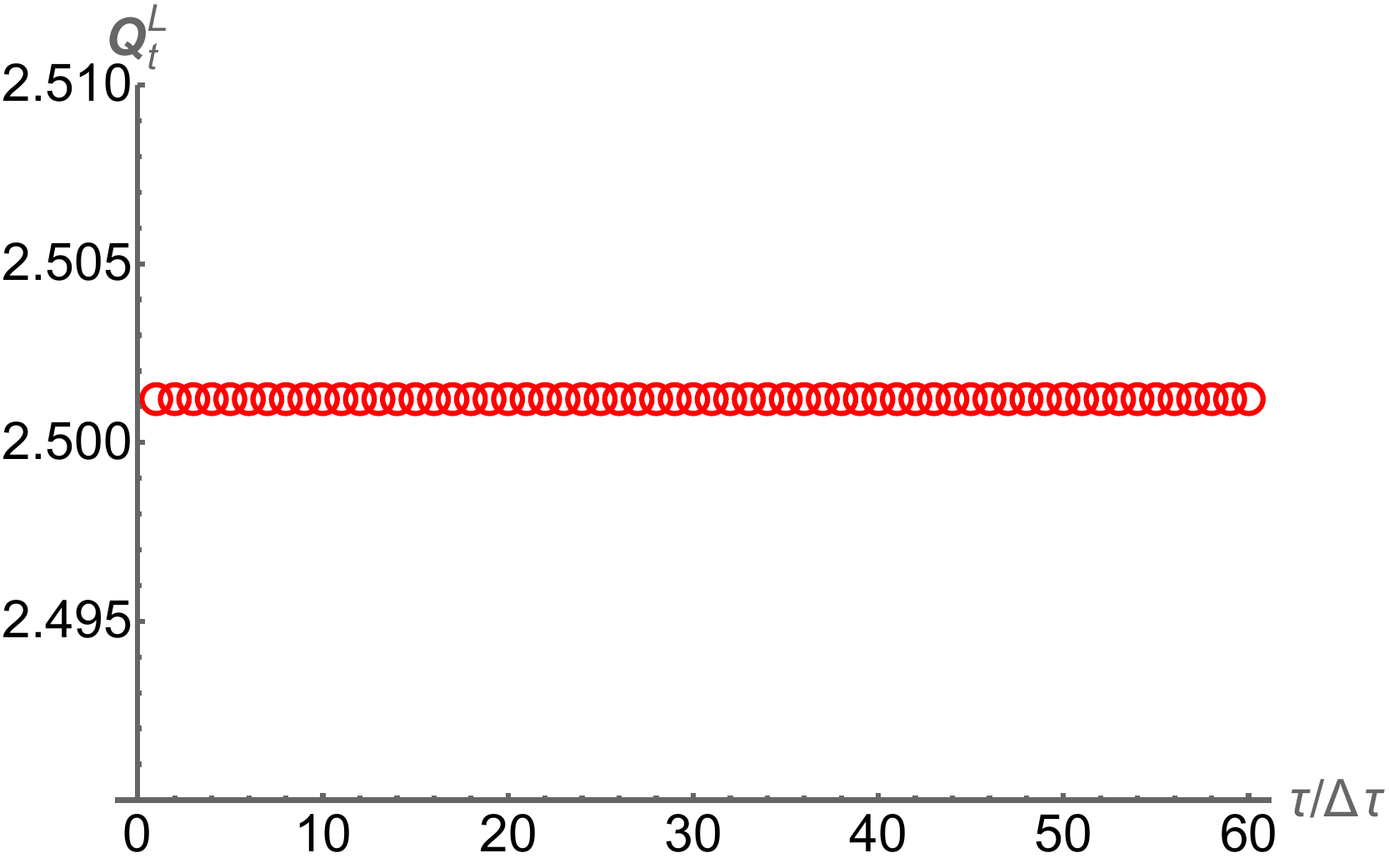}
\caption{Discrete Noether charge associated with infinitesimal time translations defined in \cref{eq:NoetherChargeDiscrLagr} and evaluated on the solution of the discretized action \cref{eq:novelactionE1p1}. Note that it is preserved exactly at each $\tau$ step. (figure from \cite{Rothkopf:2024hxi})}\label{fig:numQ}
\end{figure}

Let us continue and inspect the Noether charge of the single manifest space-time symmetry of our discretized (1+1) dimensional action, continuous time translations. In \cref{fig:numQ} we plot the values of the discrete Noether charge according to \cref{eq:NoetherChargeDiscrLagr} evaluated on the solutions shown in \cref{fig:numsol}. As one can see clearly, the values of the Noether charge are constant in $\tau$, they remain exactly conserved. If we subtract from these points the initial value of the Noether charge, we find zero within machine precision. This verifies the central claim of our approach: the retention of continuous space-time symmetries in the presence of dynamical coordinate maps after discretization, leads to an \textit{exact conservation of the associated Noether charges}.

The exact conservation of Noether charges actually provides the guiding principle for the automatic mesh refinement we observe. This can be understood already in the continuum formulation. Let us take a closer look at the expression of the Noether charge for time translations
\begin{align}
    Q_t(\tau)= \int d\sigma \Big\{ \td + \frac{1}{T}\Big( (\pp)^2\td - \pd\pp\tp \Big)\Big\} = {\rm const.}
\end{align}
In our example the solution takes on the form of two traveling wave packages. If we focus on a single rightward traveling package, its functional form reads $\phi_{\rm WP}(x-t)$. Inserting this solution into the above expression yields 
\begin{align}
    \int d\sigma \Big\{ \td + \frac{1}{T}\Big( (\pp_{\rm WP})^2(\td + \tp) \Big)\Big\} = {\rm const.}
\end{align}
Due to causality, $\td$ and $\tp$ must be positive. It follows that in order to keep the overall expression constant, the sum of $(\td + \tp)$ must adjust to the spatial gradients present in the field dynamics $(\pp_{\rm WP})^2$. In case that there are stronger gradients, the values of the temporal map derivatives must decrease. This corresponds to a finer time resolution. Conversely, if there are weaker gradients present, the temporal map derivatives must increase, which leads to an automatic coarsening of the temporal resolution. 
We have thus connected the two key observations of our proof of principle study. Introduction of dynamical coordinate maps allows to protect space-time symmetries after discretization, enforcing the conservation of associated Noether charges. In turn, this constraint on the Noether charges guides the resolution of the coordinate maps to adjust in order to capture the relevant dynamics adequately.

\section{Summary and Outlook}

In this contribution we have described how IBVPs can be formulated and solved on the level of the classical action, bypassing governing equations. The use of mimetic SBP finite difference operators gave access to Noether's theorem and associated Noether charges in a straight forward manner. In order to retain space-time symmetries, we propose to incorporate dynamical coordinate maps in addition to the conventional propagating degrees of freedom. Our proposal takes inspiration from the world-line formalism of general relativity and we derived an action for scalar wave propagation in the presence of fully dynamical coordinate maps. As a proof-of-principle we presented (1+1)d wave propagation with a reduced set of dynamical maps, incorporating only manifest time translation invariance. We confirmed the preservation of that symmetry in the discrete setting, as well as the associated conservation of the Noether charge. The coupling of field and coordinate dynamics led to a form of automatic adaptive mesh refinement.

\subsection*{Outlook}

In order to retain all space-time symmetries we need to incorporate both fully dynamic maps $t(\tau,\sigma)$ and $x(\tau,\sigma)$, which is work in progress. In the presence of fully dynamical maps, in order to find a unique solution, one must select a particular choice of coordinate representation. Analytically and in the continuum this is achieved using e.g. imposition of conformal gauge $\td\tp-\xd\xp=0$. Note that this constraint involves the derivatives of the coordinate maps and thus constitutes a non-integrable non-holonomic constraint. We have just recently finalized the action based treatment of non-holonomic constraints, based on the SKG formalism in \cite{Rothkopf:2026lru}.

We look forward to exploring the novel freedom in the treatment of spatial boundaries that arises from the presence of the new coordinate map degrees of freedom.

In the long term our goal is to extend the approach discussed in this contribution to systems with intrinsic constraints, such as Maxwell electromagnetism and to explore the conservation of Gauss' law.

\subsection*{Acknowledgements}
A.\ R.\ gladly acknowledges support by the National Research Foundation of Korea under the Emerging Researcher Grant RS-2026-25486880. A.\ R.\ thanks  Korea University for support through project K2510461 {\it Fully Dynamical Coordinate Maps for Space-Time Symmetry Preserving Lattice Field Theory} as well as project K2511131 and K2503291. J.\ N.\ was supported by Vetenskapsrådet, Sweden [award no. 2021-05484 VR] and University of Johannesburg Global
Excellence and Stature Initiative Funding. Part of this research was conducted in collaboration by WAH while visiting the Okinawa Institute of Science and Technology (OIST) through the Theoretical Sciences Visiting Program (TSVP). WAH thanks the National Research Foundation and the SA-CERN collaboration for their generous financial support during the course of this work.

\bibliographystyle{splncs04}
\bibliography{SymmetryAndMeshRefinementIBVPs}

\end{document}